\documentclass{amsart}

\theoremstyle{plain}
\newtheorem{theorem}{Theorem}[section]
\newtheorem{proposition}[theorem]{Proposition}
\newtheorem{lemma}[theorem]{Lemma}
\newtheorem{corollary}[theorem]{Corollary}
\newtheorem{example}[theorem]{Example}
\theoremstyle{definition}
\newtheorem{definition}[theorem]{Definition}
\theoremstyle{remark}
\newtheorem{remark}[theorem]{Remark}
\newtheorem*{Acknowledgments}{Acknowledgments}

\numberwithin{equation}{section}

\begin{document}

\title{Subspaces with a common complement in a Banach space}
\author{Dimosthenis Drivaliaris}
\address{Department of Financial and Management Engineering\\
University of the Aegean
31, Fostini Str.\\
82100 Chios\\
Greece}
\email{d.drivaliaris@fme.aegean.gr}
\author{Nikos Yannakakis}
\address{Department of Mathematics\\
National Technical University of Athens\\
Iroon Polytexneiou 9\\
15780 Zografou\\
Greece}
\email{nyian@math.ntua.gr}

\date{}

\begin{abstract}
We study the problem of the existence of a common algebraic complement for a pair of closed subspaces of a Banach space. We prove the following two characterizations: (1) The pairs of subspaces of a Banach space with a common complement coincide with those pairs which are isomorphic to a pair of graphs of bounded linear operators between two other Banach spaces. (2) The pairs of subspaces of a Banach space $X$ with a common complement coincide with those pairs for which there exists an involution $S$ on $X$ exchanging the two subspaces, such that $I+S$ is bounded from below on their union. Moreover we show that, in a separable Hilbert space, the only pairs of subspaces with a common complement are those which are either equivalently positioned or not completely asymptotic to one another. We also obtain characterizations for the existence of a common complement for subspaces with closed sum.
\end{abstract}
%%%%%%%%%%%%%%%%%%%%%%%%%%%%%%%%%%%%%%%%%%%%%%%%%%%%%%%%%%%%%%%%%%%%
\maketitle

\section{Introduction}
In their recent paper \cite{Lauzon} Lauzon and Treil raised the following problem: Given two closed subspaces $M$ and $N$ of a Banach space $X$, are there necessary and sufficient conditions for $M$ and $N$ to have a common algebraic complement? Recall that we say that a closed subspace $K$ is an algebraic complement (from now on just complement) of $M$ and write 
$$M\oplus K=X$$
if
$$
\begin{array}{ccc}
M\cap K=\{ 0\}&
\mathrm{ and }&
M+K=X.
\end{array}
$$
So what one is looking for are conditions equivalent to the existence of a third closed subspace $K$ of $X$, which we will call a common complement of $M$ and $N$ in $X$, with
$$M\oplus K=N\oplus K=X.$$
%%%%%%%%%%%%%%%%%%%%%%%%%%%%%%%%%%%%%%%%%%%%%%%%%%%%%%%%%%%%%%%%%%%%

It is well known that for a finite dimensional Banach space $X$ such a subspace $K$ exists if and only if $M$ and $N$ have equal dimensions. If we move to an infinite dimensional Hilbert space things are much more complicated and equality of dimensions and codimensions is necessary (in general in a Banach space subspaces with a common complement are isomorphic), but no longer sufficient; an easy way to see this is to let $M$ and $N$ have both infinite dimensions and codimensions and $N$ be a proper subspace of $M$. This should come as no surprise, since the equality of dimensions and codimensions is necessary and sufficient only if we are looking for something less, namely not for an algebraic but for a topological common complement. Recall that $K$ is a topological complement of $M$ in $X$ if 
$$
\begin{array}{ccc}
M\cap K=\{ 0\}&
\mathrm{ and }&
\overline{M+K}=X.
\end{array}
$$
If $M$ and $N$ have a common topological complement $K$ and $P_M$, $P_N$ and $P_K$ are the orthogonal projections on $M$, $N$ and $K$ respectively, then
$$
\begin{array}{ccc}
P_M\wedge P_K=0=P_N\wedge P_K&
\mathrm{and}&
P_M\vee P_K=I=P_N\vee P_K.
\end{array}
$$
Pairs of projections with the property just described for $P_M$ and $P_N$ are called perspective. The study of perspective projections goes back to Kaplansky \cite[Theorem 6.6]{Kaplansky}. Fillmore proved \cite[Theorem 1]{Fillmore} that two orthogonal projections are perspective in the projection lattice of a von Neumann algebra if and only if they are unitarily equivalent. A spatial interpretation of Fillmore's result gives us what we stated above: $M$ and $N$ have a common topological complement if and only if they have equal dimensions and codimensions. For further results on perspective projections we refer the reader to \cite{Elliott,Holland,Maeda}. 
%%%%%%%%%%%%%%%%%%%%%%%%%%%%%%%%%%%%%%%%%%%%%%%%%%%%%%%%%%%%%%%%%%%%

A stronger, than the equality of dimensions and codimensions, condition which could be a possible candidate for the characterization of the existence of a common complement is the following
\begin{equation}
\label{dimensions}
\mathrm{dim} (M\ominus (M\cap N))=\mathrm{dim} (N\ominus (M\cap N)).
\end{equation}
It turns out, see the example in \cite[Section 4]{Lauzon}, that (\ref{dimensions}) is also not sufficient.
%%%%%%%%%%%%%%%%%%%%%%%%%%%%%%%%%%%%%%%%%%%%%%%%%%%%%%%%%%%%%%%%%%%%

In \cite[Theorem 0.1]{Lauzon} Lauzon and Treil obtained the following characterization:
\begin{theorem}
\label{Treil}
Let $M$ and $N$ be two closed subspaces of a Hilbert space $X$ and $G$ be the restriction on $M$ of the orthogonal projection $P_N$. Then $M$ and $N$ have a common complement if and only if
\begin{equation}
\label{Lauzon}
\mathrm{dim}(M\cap N^\perp)+\mathrm{dim}(\mathcal{E}((0,1-\varepsilon))(M))= \mathrm{dim}(M^\perp\cap N)+\mathrm{dim}(\mathcal{E}((0,1-\varepsilon))(M)),
\end{equation}
for some $\varepsilon >0$ (for all sufficiently small $\varepsilon >0$), where $\mathcal{E}(\cdot)$ is the spectral measure of the operator $G^*G$.
\end{theorem}
A different proof of this result for a separable Hilbert space $X$ was given by Du and Deng in \cite{Du}. Related results, also for the separable case, using the minimal number of segments required to connect two homotopic projections, can be found in \cite[Proposition 6.2 and Theorem 6.3]{Giol}. 
%%%%%%%%%%%%%%%%%%%%%%%%%%%%%%%%%%%%%%%%%%%%%%%%%%%%%%%%%%%%%%%%%%%%

Note that the above theorem implies that if two subspaces $M$ and $N$ of a Hilbert space are, in the terminology of \cite{Davis}, equivalently positioned, i.e.\ if
\begin{equation}
\label{equiv}
\mathrm{dim}(M\cap N^\perp)= \mathrm{dim}(M^\perp\cap N),
\end{equation}
then they have a common complement (as we will see in Example \ref{example} the converse is not true). 
This trivial observation allows us to settle the case of finite dimensional and codimensional subspaces. If $M$ and $N$ are finite dimensional, then (\ref{equiv}) holds if and only if $\mathrm{dim}(M)= \mathrm{dim}(N)$ \cite[p.\ 175]{Davis}. Moreover $M^\perp$ and $N^\perp$ are equivalently positioned if and only if $M$ and $N$ are. Hence for finite dimensional or codimensional subspaces of a Hilbert space the existence of a common complement is equivalent to their dimensions or codimensions being equal. We should add here that equivalently positioned subspaces behave so well because there exists a symmetry exchanging them, see \cite[p. 389]{Dixmier}, and a common complement is the orthogonal complement of its axis. This observation will play an important role in Section \ref{involutions}.
%%%%%%%%%%%%%%%%%%%%%%%%%%%%%%%%%%%%%%%%%%%%%%%%%%%%%%%%%%%%%%%%%%%%

If $M$ and $N$ have a common complement, but are not equivalently positioned, then condition (\ref{Lauzon}) implies that the dimension of $\mathcal{E}((0,1-\varepsilon))(M)$ is infinite. Using this observation we shall see that, at least for the separable case, the existence of a common complement for two subspaces has to do either with being in a ``good relative position'' or with both of them having infinite dimensional closed subspaces ``away'' from the other. The latter is precisely the definition of subspaces not completely asymptotic to one another (see Definition \ref{compasymp}). Hence we can say that among the ``not so nicely positioned'' subspaces (the not equivalently positioned ones) the only pairs that have a common complement, always in the separable case, are those which are not completely asymptotic to one another. Concluding, it is interesting to note that in Section \ref{graphs} we will see that up to isomorphism all subspaces with a common complement in a Hilbert space are equivalently positioned.
%%%%%%%%%%%%%%%%%%%%%%%%%%%%%%%%%%%%%%%%%%%%%%%%%%%%%%%%%%%%%%%%%%%%

Another characterization of subspaces of a Hilbert space with a common complement is the following by Longstaff and Panaia \cite[Proposition 1]{Longstaff}: 
\begin{theorem}
\label{Longstaff}
Let $M$ and $N$ be subspaces of a Hilbert space $X$ with\\$\overline{M+N}=X$ and $M\cap N=\{ 0\}$. Then the following are equivalent:
\begin{enumerate}
\item $M$ and $N$ have a common complement in $X$.
\item The pair $\left\{ M,N\right\}$ is similar to a pair of the form $\left\{ \mathrm{Gr}(T),\mathrm{Gr}(S)\right\}$, where $T$ and $S$ are operators on some Hilbert space $X_1$.
\item  The pair $\left\{ M,N\right\}$ is similar to a pair of subspaces in generic position (see Section 4 for the definition).
\end{enumerate}
\end{theorem}
Note that this characterization applies only to topologically complementary subspaces and that, whereas the characterization of Lauzon and Treil concentrates on dimensions, this one focuses on the relative position of the pair. 
%%%%%%%%%%%%%%%%%%%%%%%%%%%%%%%%%%%%%%%%%%%%%%%%%%%%%%%%%%%%%%%%%%%%

Sufficient conditions for the existence of a common complement in a Banach space, using notions of distance between subspaces, can be found in the papers of Berkson \cite[Theorem 5.2 and Corollary 5.7]{Berkson} and of Dirr, Rako\v{c}evi\'{c} and Wimmer \cite[Theorem 3.1]{Dirr}. 
%%%%%%%%%%%%%%%%%%%%%%%%%%%%%%%%%%%%%%%%%%%%%%%%%%%%%%%%%%%%%%%%%%%%

Our main aim in this paper is to obtain characterizations for pairs (or families) of subspaces with a common complement in a Banach space. Moreover we shall examine their implications in the Hilbert space setting and get
more information about the existence of a common complement there. We should note that the tools used by Lauzon and Treil are not suitable for this more general setting. On the other hand the characterization of Longstaff and Panaia and the discussion about equivalently positioned subspaces indicate that it might be fruitful to treat our problem as one of relative position of a pair of subspaces in a Banach space. We remind the reader that the study of the relative position of subspaces of a Hilbert space was initiated by Dixmier in \cite{Dixmier} and has been going on ever since (see for example the papers of Davis \cite{Davis}, Araki \cite{Araki}, Halmos \cite{Halmos}, Davis and Kahan \cite{Davis1}, Longstaff and Panaia \cite{Longstaff} and Enomoto and Watatani \cite{Enomoto}). 
%%%%%%%%%%%%%%%%%%%%%%%%%%%%%%%%%%%%%%%%%%%%%%%%%%%%%%%%%%%%%%%%%%%%

Our paper is organized as follows: In the beginning of Section \ref{preliminaries} we present some simple facts about subspaces with a common complement. We believe that this discussion shows some of the problems that may occur when we study the existence of a common complement in the Banach space setting. As we have already said, if two subspaces have a common complement, then they are isomorphic, whereas the converse is not true. We conclude this section with a condition for the isomorphism taking the one subspace to the other which characterizes subspaces with a common complement among the isomorphic ones. 
%%%%%%%%%%%%%%%%%%%%%%%%%%%%%%%%%%%%%%%%%%%%%%%%%%%%%%%%%%%%%%%%%%%%

In Section \ref{closed} we discuss the relation between the existence of a common complement for two subspaces and their sum being closed or them being completely asymptotic. To this end we first deal with the case of a separable Hilbert space. Using the results of Lauzon and Treil we show that $M$ and $N$ have a common complement if and only if either they are equivalently positioned or they are not completely asymptotic to one another. An immediate consequence is that, in a separable Hilbert space, two subspaces $M$ and $N$ with closed sum have a common complement if and only if the equality of dimensions (\ref{dimensions}) holds. Applying our results from Section \ref{preliminaries}, we extend this last result to arbitrary Banach spaces. 
%%%%%%%%%%%%%%%%%%%%%%%%%%%%%%%%%%%%%%%%%%%%%%%%%%%%%%%%%%%%%%%%%%%%

In Section \ref{graphs} we extend characterizations of pairs of subspaces in various relative positions as those isomorphic to pairs of graphs of operators with certain properties, due to Araki \cite{Araki}, Halmos \cite{Halmos}, Papadakis \cite{Papadakis} and Longstaff and Panaia \cite{Longstaff}, to pairs of subspaces with a common complement. We characterize pairs (or families) of closed subspaces of a Banach space $X$ that have a common complement as those which are isomorphic to pairs (or families) of graphs of bounded linear operators between two Banach spaces $X_1$ and $X_2$. An interesting consequence of this is that up to isomorphism all pairs of subspaces of a Hilbert space with a common complement are equivalently positioned. We conclude by characterizing pairs $\{ M, N\}$ of subspaces of a Hilbert space for which $M^{\perp}$ is a common complement.
%%%%%%%%%%%%%%%%%%%%%%%%%%%%%%%%%%%%%%%%%%%%%%%%%%%%%%%%%%%%%%%%%%%%

In the final section we extend results of Dixmier \cite{Dixmier} and Davis \cite{Davis}, characterizing pairs of subspaces in various relative positions via the existence of certain kinds of symmetries interchanging them, to pairs of subspaces with a common complement. We show that for two subspaces $M$ and $N$ of a Banach space $X$ with $\overline{M+N}=X$, having a common complement is equivalent to the existence of an involution $S$ on $X$ (i.e.\ a bounded linear operator with $S^2=I$) which exchanges $M$ and $N$, such that $I+S$ is bounded from below on their union. 
%%%%%%%%%%%%%%%%%%%%%%%%%%%%%%%%%%%%%%%%%%%%%%%%%%%%%%%%%%%%%%%%%%%%

\section{Preliminaries and First Results}
\label{preliminaries}
%%%%%%%%%%%%%%%%%%%%%%%%%%%%%%%%%%%%%%%%%%%%%%%%%%%%%%%%%%%%%%%%%%%%

We start by clarifying some notational matters. All the results of the paper are true both for real and for complex Banach and Hilbert spaces. Throughout the paper all subspaces are considered to be closed. We will call an operator $U:X_1\rightarrow X_2$ an isomorphism if it is an injective bounded linear operator with closed range.
%%%%%%%%%%%%%%%%%%%%%%%%%%%%%%%%%%%%%%%%%%%%%%%%%%%%%%%%%%%%%%%%%%%%

Let $M$ and $N$ be subspaces of a Banach space. We will denote the pair of those two subspaces with no order specified by $\left\{ M,N\right\}$. If $X_1$ and $X_2$ are Banach spaces and $\left\{ M_1,N_1\right\}$ and $\left\{ M_2,N_2\right\}$ are pairs of subspaces of $X_1$ and $X_2$ respectively we will say that the pairs $\left\{ M_1,N_1\right\}$ and $\left\{ M_2,N_2\right\}$ are isomorphic/isometrically isomorphic if there exists an isomorphism onto/isometric isomorphism onto taking the one pair to the other. If $X_1$ and $X_2$ are Hilbert spaces we will use the expressions similar and unitarily equivalent instead of isomorphic and isometrically isomorphic. 
%%%%%%%%%%%%%%%%%%%%%%%%%%%%%%%%%%%%%%%%%%%%%%%%%%%%%%%%%%%%%%%%%%%%

If $M$ and $K$ are complementary we will denote the projection onto $M$ parallel to $K$ by $P_{M\| K}$. If $M$ is a subspace of a Hilbert space $X$ we will denote the orthogonal projection onto $M$ by $P_M$. Moreover if $X$ is a Hilbert space, $M$ is a subspace of $X$ and $L$ is a subspace of $M$, then we will denote the orthogonal complement of $L$ in $M$ by $M\ominus L$.
%%%%%%%%%%%%%%%%%%%%%%%%%%%%%%%%%%%%%%%%%%%%%%%%%%%%%%%%%%%%%%%%%%%%

We describe now some simple properties of subspaces with a common complement.
%%%%%%%%%%%%%%%%%%%%%%%%%%%%%%%%%%%%%%%%%%%%%%%%%%%%%%%%%%%%%%%%%%%%

\begin{proposition}
\label{proposition1}
Let $X$ be a Banach space, $M$ and $N$ be subspaces of $X$ and $Y$ be a complemented subspace of $X$ containing both $M$ and $N$. Then the following are equivalent:
\begin{enumerate}
\item $M$ and $N$ have a common complement in $X$.
\item $M$ and $N$ have a common complement in $Y$.
\end{enumerate}
\end{proposition}
%%%%%%%%%%%%%%%%%%%%%%%%%%%%%%%%%%%%%%%%%%%%%%%%%%%%%%%%%%%%%%%%%%%%

\begin{proof}
If $K$ is a common complement of $M$ and $N$ in $X$, then it is easy to see that $Y\cap K$ is a common complement of $M$ and $N$ in $Y$. On the other hand if $K$ is a common complement of $M$ and $N$ in $Y$ and $Z$ is a complement of $Y$ in $X$, then $K\oplus Z$ is a common complement of $M$ and $N$ in $X$. 
\end{proof}
%%%%%%%%%%%%%%%%%%%%%%%%%%%%%%%%%%%%%%%%%%%%%%%%%%%%%%%%%%%%%%%%%%%%

\begin{remark}
\label{remark1}
The direction (1)$\Rightarrow$(2) is true even if $Y$ is not complemented in $X$, whereas the direction (2)$\Rightarrow$(1) is not in general true if $Y$ is not complemented in $X$, since $M$ and $N$ may be non-complemented subspaces of $X$ and have a common complement in some non-complemented subspace $Y$ containing both of them.
\end{remark}
%%%%%%%%%%%%%%%%%%%%%%%%%%%%%%%%%%%%%%%%%%%%%%%%%%%%%%%%%%%%%%%%%%%%

A straightforward corollary of Proposition \ref{proposition1} is that if $M$ and $N$ are finite dimensional subspaces of a Banach space $X$, then $M$ and $N$ have a common complement in $X$ if and only if $\mathrm{dim}(M)=\mathrm{dim}(N)$ (just take $Y=M+N$ and apply the result for the finite dimensional case).
%%%%%%%%%%%%%%%%%%%%%%%%%%%%%%%%%%%%%%%%%%%%%%%%%%%%%%%%%%%%%%%%%%%%

We will later see that some results for the existence of a common complement hold if we assume that the subspaces have dense sum. In those cases the following corollary, which is a Banach space version of \cite[Proposition 1.5]{Lauzon}, is very useful.   
%%%%%%%%%%%%%%%%%%%%%%%%%%%%%%%%%%%%%%%%%%%%%%%%%%%%%%%%%%%%%%%%%%%%

\begin{corollary}
\label{corollary1}
Let $X$ be a Banach space and $M$ and $N$ be subspaces of $X$ with $\overline{M+N}$ complemented in $X$. Then $M$ and $N$ have a common complement in $X$ if and only if they have a common complement in $\overline{M+N}$.
\end{corollary}
%%%%%%%%%%%%%%%%%%%%%%%%%%%%%%%%%%%%%%%%%%%%%%%%%%%%%%%%%%%%%%%%%%%

In the following two propositions we describe what happens with the existence of a common complement if we ``remove'' a common part from $M$ and $N$. 
%%%%%%%%%%%%%%%%%%%%%%%%%%%%%%%%%%%%%%%%%%%%%%%%%%%%%%%%%%%%%%%%%%%%

\begin{proposition}
\label{proposition3}
Let $X$ be a Banach space, $M$ and $N$ be subspaces of $X$ and $L$ be a subspace both of $M$ and of $N$ which is complemented in $X$. If $P:X\rightarrow X$ is a bounded projection with $P(X)=L$, then the following are equivalent:
\begin{enumerate}
\item $M$ and $N$ have a common complement in $X$.
\item $(I-P)(M)$ and $(I-P)(N)$ have a common complement in $X$.
\item $(I-P)(M)$ and $(I-P)(N)$ have a common complement in $(I-P)(X)$.
\end{enumerate}
\end{proposition}
%%%%%%%%%%%%%%%%%%%%%%%%%%%%%%%%%%%%%%%%%%%%%%%%%%%%%%%%%%%%%%%%%%%%

\begin{proof}
If $K$ is a common complement of $M$ and $N$ in $X$, then $L\oplus K$ is a common complement of $(I-P)(M)$ and $(I-P)(N)$ in $X$ and so we get (1)$\Rightarrow$(2). (2)$\Rightarrow$(3) follows from Proposition \ref{proposition1}. For (3)$\Rightarrow$(1) observe that if $K$ is a common complement of $(I-P)(M)$ and $(I-P)(N)$ in $(I-P)(X)$, then $K$ is also a common complement of $M$ and $N$ in $X$. 
\end{proof}
%%%%%%%%%%%%%%%%%%%%%%%%%%%%%%%%%%%%%%%%%%%%%%%%%%%%%%%%%%%%%%%%%%%%

If the subspace $L$ is just complemented in $M$ and $N$, then we get the following.
%%%%%%%%%%%%%%%%%%%%%%%%%%%%%%%%%%%%%%%%%%%%%%%%%%%%%%%%%%%%%%%%%%%%

\begin{proposition}
\label{corollary4}
Let $X$ be a Banach space and $M,\,N$ and $L$ be subspaces of $X$ such that $L$ is complemented both in $M$ and in $N$. If $M=L\oplus M_1$ and $N=L\oplus N_1$, then the following are equivalent:
\begin{enumerate}
\item $M$ and $N$ have a common complement in $X$.
\item $M_1$ and $N_1$ have a common complement in $X$ in which $L$ is complemented.
\end{enumerate}
\end{proposition}
%%%%%%%%%%%%%%%%%%%%%%%%%%%%%%%%%%%%%%%%%%%%%%%%%%%%%%%%%%%%%%%%%%%%

\begin{proof}
If $K$ is a common complement of $M$ and $N$ in $X$, then $L\oplus K$ is a common complement of $M_1$ and $N_1$ in $X$ in which $L$ is complemented. On the other hand if $K=L\oplus K'$ is a common complement of $M_1$ and $N_1$ in $X$ in which $L$ is complemented, then $K'$ is a common complement of $M$ and $N$.
\end{proof}
%%%%%%%%%%%%%%%%%%%%%%%%%%%%%%%%%%%%%%%%%%%%%%%%%%%%%%%%%%%%%%%%%%%%

\begin{remark}
\label{remark1a}
The direction (2)$\Rightarrow$(1) is not in general true if we omit ``in which $L$ is complemented'' in (2), since in that case $L$ may be a non-complemented subspace of $X$.
\end{remark}
%%%%%%%%%%%%%%%%%%%%%%%%%%%%%%%%%%%%%%%%%%%%%%%%%%%%%%%%%%%%%%%%%%%%

In general it is not true that if $M$ and $N$ have a common complement, then their intersection is a complemented subspace of $X$. If it is we get the following immediate corollary of Proposition \ref{proposition3}.
%%%%%%%%%%%%%%%%%%%%%%%%%%%%%%%%%%%%%%%%%%%%%%%%%%%%%%%%%%%%%%%%%%%%

\begin{corollary}
\label{corollary3}
Let $X$ be a Banach space and $M$ and $N$ be subspaces of $X$ such that $M\cap N$ is complemented in $X$. If $P:X\rightarrow X$ is a bounded projection with $P(X)=M\cap N$, then $M$ and $N$ have a common complement in $X$ if and only if $(I-P)(M)$ and $(I-P)(N)$ have a common complement in $(I-P)(X)$.
\end{corollary}
%%%%%%%%%%%%%%%%%%%%%%%%%%%%%%%%%%%%%%%%%%%%%%%%%%%%%%%%%%%%%%%%%%%%

Combining Corollaries \ref{corollary1} and \ref{corollary3} we get that if $M$ and $N$ are subspaces of a Banach space $X$ such that $\overline{M+N}$ and $M\cap N$ are complemented in $X$ and $P$ is a bounded projection onto $M\cap N$, then $M$ and $N$ have a common complement in $X$ if and only if $M_1=(I-P)(M)$ and $N_1=(I-P)(N)$ have a common complement in $(I-P)(\overline{M+N})=\overline{M_1+N_1}$. This is quite useful since we end up with two subspaces with dense sum and trivial intersection. In particular when we study the existence of a common complement in a Hilbert space we may assume that we are talking about subspaces with dense sum and trivial intersection.
%%%%%%%%%%%%%%%%%%%%%%%%%%%%%%%%%%%%%%%%%%%%%%%%%%%%%%%%%%%%%%%%%%%%

It is straightforward to see that the existence of a common complement is preserved under isomorphisms.
%%%%%%%%%%%%%%%%%%%%%%%%%%%%%%%%%%%%%%%%%%%%%%%%%%%%%%%%%%%%%%%%%%%%

\begin{proposition}
\label{proposition5}
Let $X_1$ and $X_2$ be Banach spaces, $M$ and $N$ be subspaces of $X_1$ and $U:X_1\rightarrow X_2$ be an isomorphism. Then $M$ and $N$ have a common complement in $X_1$ if and only if $U(M)$ and $U(N)$ have a common complement in $U(X_1)$.
\end{proposition}
%%%%%%%%%%%%%%%%%%%%%%%%%%%%%%%%%%%%%%%%%%%%%%%%%%%%%%%%%%%%%%%%%%%%

Moreover the existence of a common complement is also preserved under taking annihilators.
%%%%%%%%%%%%%%%%%%%%%%%%%%%%%%%%%%%%%%%%%%%%%%%%%%%%%%%%%%%%%%%%%%%%

\begin{proposition}
\label{proposition6}
Let $X$ be a Banach space and $M$ and $N$ be subspaces of $X$. 
\begin{enumerate}
\item If $M$ and $N$ have a common complement in $X$, then $M^{\perp}$ and $N^{\perp}$ have a common complement in $X^*$.
\item If $X$ is reflexive and $M^{\perp}$ and $N^{\perp}$ have a common complement in $X^*$, then $M$ and $N$ have a common complement in $X$.
\end{enumerate}
\end{proposition}
%%%%%%%%%%%%%%%%%%%%%%%%%%%%%%%%%%%%%%%%%%%%%%%%%%%%%%%%%%%%%%%%%%%%

We should note that (2) is not in general true if $X$ is not reflexive since in this case $M$ and $N$ may not even be complemented.
%%%%%%%%%%%%%%%%%%%%%%%%%%%%%%%%%%%%%%%%%%%%%%%%%%%%%%%%%%%%%%%%%%%%

Combining the previous proposition with what we said after Remark \ref{remark1} we get that if $X$ is a reflexive Banach space and $M$ and $N$ are finite codimensional subspaces of $X$, then $M$ and $N$ have a common complement in $X$ if and only if $\mathrm{codim}(M)=\mathrm{codim}(N)$.
%%%%%%%%%%%%%%%%%%%%%%%%%%%%%%%%%%%%%%%%%%%%%%%%%%%%%%%%%%%%%%%%%%%%

The existence of a common topological complement is a transitive relation \cite[Corollary, p.\ 386]{Fillmore}. On the other hand the relation of being equivalently positioned is not transitive \cite[p.\ 175]{Davis}. The example from \cite{Davis} can also be used to show that the existence of a common complement is not transitive. Let $X=l^2(\mathbb{Z})$, $\{ e_n\,|\,n\in\mathbb{Z}\}$ be the standard orthonormal basis of $X$,
$$
\begin{array}{ccc}
M=\overline{\mathrm{span}}\{ e_n\,|\,n\geq 0\},&
N=\overline{\mathrm{span}}\{ e_n\,|\,n\leq -1\},&
L=\overline{\mathrm{span}}\{ e_n\,|\,n\geq 1\}.
\end{array}
$$
Then  
$$
\begin{array}{ccc}
\mathrm{dim}(M\cap N^\perp)= \mathrm{dim}(M^\perp\cap N)&
\mathrm{and}&
\mathrm{dim}(N\cap L^\perp)= \mathrm{dim}(N^\perp\cap L)
\end{array}
$$
and thus $M$ and $N$ have a common complement and $N$ and $L$ have a common complement. On the other hand, since $L$ is a proper subspace of $M$, $M$ and $L$ do not have a common complement. 
%%%%%%%%%%%%%%%%%%%%%%%%%%%%%%%%%%%%%%%%%%%%%%%%%%%%%%%%%%%%%%%%%%%%

We move now to the next part of this section. Lauzon and Treil proved \cite[Proposition 1.3]{Lauzon} that if $M$ and $N$ are subspaces of a Banach space $X$, $M$ and $N$ have a common complement in $X$ if and only if there exists a bounded projection $P:X\rightarrow X$ with $P(X)=N$ such that $P|_M:M\rightarrow N$ is an isomorphism onto. It is obvious that this leads to the following corollary.
%%%%%%%%%%%%%%%%%%%%%%%%%%%%%%%%%%%%%%%%%%%%%%%%%%%%%%%%%%%%%%%%%%%%

\begin{corollary}
\label{corollary6}
If $M$ and $N$ have a common complement in $X$, then $M$ and $N$ are isomorphic.
\end{corollary}
%%%%%%%%%%%%%%%%%%%%%%%%%%%%%%%%%%%%%%%%%%%%%%%%%%%%%%%%%%%%%%%%%%%%

Combining Corollary \ref{corollary6} with the simple properties of subspaces with a common complement we get that if two subspaces have a common complement, then other subspaces related to them are also isomorphic (for example, if $M$ and $N$ have a common complement, then their annihilators are isomorphic). In particular combining Corollaries \ref{corollary3} and \ref{corollary6} we get the following:
%%%%%%%%%%%%%%%%%%%%%%%%%%%%%%%%%%%%%%%%%%%%%%%%%%%%%%%%%%%%%%%%%%%%

\begin{corollary}
\label{corollary10}
Let $X$ be a Banach space and $M$ and $N$ be subspaces of $X$ such that $M\cap N$ is complemented in $X$ and let $P:X\rightarrow X$ be a bounded projection with $P(X)=M\cap N$. If $M$ and $N$ have a common complement in $X$, then $(I-P)(M)$ and $(I-P)(N)$ are isomorphic. 
\end{corollary}
%%%%%%%%%%%%%%%%%%%%%%%%%%%%%%%%%%%%%%%%%%%%%%%%%%%%%%%%%%%%%%%%%%%%

For Hilbert spaces, what we just said implies, by taking $P=P_{M\cap N}$, the result of \cite[Corollary 1.4]{Lauzon}: Let $X$ be a Hilbert space and $M$ and $N$ be subspaces of $X$. If $M$ and $N$ have a common complement in $X$, then
$$\mathrm{dim} (M\ominus (M\cap N))=\mathrm{dim} (N\ominus (M\cap N)).$$
%%%%%%%%%%%%%%%%%%%%%%%%%%%%%%%%%%%%%%%%%%%%%%%%%%%%%%%%%%%%%%%%%%%%

As we mentioned in the Introduction, the converse of Corollaries \ref{corollary6} and \ref{corollary10} is not true. 
%%%%%%%%%%%%%%%%%%%%%%%%%%%%%%%%%%%%%%%%%%%%%%%%%%%%%%%%%%%%%%%%%%%%

In the following proposition we give a condition on the isomorphism between two isomorphic subspaces of a Banach space equivalent to them having a common complement in the closure of their sum. This result will be used in the following section where we discuss the existence of a common complement for subspaces with closed sum.
%%%%%%%%%%%%%%%%%%%%%%%%%%%%%%%%%%%%%%%%%%%%%%%%%%%%%%%%%%%%%%%%%%%%

\begin{proposition}
\label{proposition10}
Let $X$ be a Banach space and $M$ and $N$ be subspaces of $X$. Then the following are equivalent:
\begin{enumerate}
\item $M$ and $N$ have a common complement in $\overline{M+N}$.
\item There exists an isomorphism onto $U:M\rightarrow N$, with $U|_{M\cap N}=I_{M\cap N}$, and $C\geq 1$ such that 
$$\| Ux+y\|\leq C\| x+y\| ,$$
for all $x\in M$ and $y\in N$.
\end{enumerate}
\end{proposition}
%%%%%%%%%%%%%%%%%%%%%%%%%%%%%%%%%%%%%%%%%%%%%%%%%%%%%%%%%%%%%%%%%%%%

\begin{proof}
(1)$\Rightarrow$(2): Since $M$ and $N$ have a common complement in $\overline{M+N}$, by \cite[Proposition 1.3]{Lauzon}, there exists a bounded projection
$$P:\overline{M+N}\rightarrow \overline{M+N}$$
with $P(\overline{M+N})=N$ such that $P|_M:M\rightarrow N$ is an isomorphism onto. Let $U=P|_M$. If $x\in M\cap N$, then obviously $Ux=P|_Mx=x$ and thus $U|_{M\cap N}=I_{M\cap N}$. Moreover if $x\in M$ and $y\in N$, then
$$
\|Ux+y\|=\|P|_M x+y\|=\|P(x+y)\|\leq \|P\|\;\|x+y\| .
$$
Hence, for $C=\| P\|\geq 1$, we have
$$\| Ux+y\|\leq C\| x+y\| ,$$
for all $x\in M, y\in N$.

\noindent (2)$\Rightarrow$(1): Define
$$P:M+N\rightarrow N$$
with
$$P(x+y)=Ux+y,$$
for all $x\in M$ and $y\in N$. To see that $P$ is well-defined let $x_1,\, x_2\in M$ and $y_1,\, y_2\in N$ with $x_1+y_1=x_2+y_2$. Then
$$
x_1-x_2=y_2-y_1\in M\cap N.
$$
Therefore, since $U|_{M\cap N}=I_{M\cap N}$,
$$
Ux_1-Ux_2=y_2-y_1
$$
and hence
$$
P(x_1+y_1)= Ux_1+y_1= Ux_2+y_2= P(x_2+y_2).
$$
Since $U$ is an isomorphism onto $N$, $P$ is a projection onto $N$.
Moreover if $x\in M$ and $y\in N$, then
$$
\| P(x+y)\|=\| Ux+y\|\leq C\| x+y\|
$$
and thus $P$ is bounded. Extending $P$ by continuity to the whole of $\overline{M+N}$ we get a bounded projection $\hat{P}:\overline{M+N}\rightarrow\overline{M+N}$ with $\hat{P}(\overline{M+N})=N$. Since
$$\hat{P}|_M=P|_M=U,$$
$\hat{P}|_M$ is an isomorphism. Therefore, by \cite[Proposition 1.3]{Lauzon}, $M$ and $N$ 
have a common complement in $\overline{M+N}$.
\end{proof}
%%%%%%%%%%%%%%%%%%%%%%%%%%%%%%%%%%%%%%%%%%%%%%%%%%%%%%%%%%%%%%%%%%%%

Note that the existence of an isomorphism like the one described above does not in general imply the existence of a common complement of $M$ and $N$ in the whole of $X$. 
%%%%%%%%%%%%%%%%%%%%%%%%%%%%%%%%%%%%%%%%%%%%%%%%%%%%%%%%%%%%%%%%%%%%

\section{Subspaces with closed sum}
\label{closed}
%%%%%%%%%%%%%%%%%%%%%%%%%%%%%%%%%%%%%%%%%%%%%%%%%%%%%%%%%%%%%%%%%%%%

One of the main complications in infinite dimensions is that the sum of two subspaces may not be closed. In general, subspaces with a common complement need not have a closed sum. To see that let $X=l^2(\mathbb{N})$, $\{ e_n\,|\,n\in\mathbb{N}\}$ be the standard orthonormal basis of $X$,
$$
\begin{array}{ccc}
M=\overline{\mathrm{span}}\left\{e_{2n-1}\,|\,n\geq 1\right\}&
\mathrm{and}&
N=\overline{\mathrm{span}}\left\{\left.\sqrt{\frac{n-1}{n}}\,e_{2n-1}+\sqrt{\frac{1}{n}}\,e_{2n}\,\right|
\,n\geq 1\right\}.
\end{array}
$$ 
Then 
$$M\cap N^{\perp}=\{ 0\}=M^{\perp}\cap N$$
and thus, by Theorem \ref{Treil}, $M$ and $N$ have a common complement, but $M+N$ is not closed. In this section we will discuss the connections between the existence of a common complement for two subspaces and their sum being closed.
%%%%%%%%%%%%%%%%%%%%%%%%%%%%%%%%%%%%%%%%%%%%%%%%%%%%%%%%%%%%%%%%%%%%

We start with a theorem that shows that in a separable Hilbert space two subspaces have a common complement if and only if they are either equivalently positioned or ``do not behave very bad with respect to their sum being closed''. To make the last phrase precise we need the notion of completely asymptotic subspaces which is due to Dixmier \cite[p. 23]{Dixmier2}.
%%%%%%%%%%%%%%%%%%%%%%%%%%%%%%%%%%%%%%%%%%%%%%%%%%%%%%%%%%%%%%%%%%%%

\begin{definition}
\label{compasymp}
Let $X$ be a separable Hilbert space and $M$ and $N$ be subspaces of $X$. We will say that $M$ is completely asymptotic to $N$ if for any infinite dimensional subspace $M_1$ of $M$ which is disjoint from $N$, $M_1+N$ is not closed. Moreover we will say that $M$ and $N$ are completely asymptotic if they are completely asymptotic to one another.
\end{definition}
%%%%%%%%%%%%%%%%%%%%%%%%%%%%%%%%%%%%%%%%%%%%%%%%%%%%%%%%%%%%%%%%%%%%%%

\begin{remark}
\label{asym}
By \cite[Theorem 1.4]{Dixmier2} two subspaces $M$ and $N$ are not completely asymptotic to one another if and only if both $(I-P_M)(N)$ and $(I-P_N)(M)$ contain an infinite dimensional subspace.
\end{remark}
%%%%%%%%%%%%%%%%%%%%%%%%%%%%%%%%%%%%%%%%%%%%%%%%%%%%%%%%%%%%%%%%%%%%%

We can now state the theorem that we described above.
%%%%%%%%%%%%%%%%%%%%%%%%%%%%%%%%%%%%%%%%%%%%%%%%%%%%%%%%%%%%%%%%%%%%%

\begin{theorem}
\label{asymptotic}
Let $X$ be a separable Hilbert space and $M$ and $N$ be subspaces of $X$. Then $M$ and $N$ have a common complement if and only if they are either equivalently positioned or not completely asymptotic to one another. 
\end{theorem}
%%%%%%%%%%%%%%%%%%%%%%%%%%%%%%%%%%%%%%%%%%%%%%%%%%%%%%%%%%%%%%%%%%%%%

For the proof of the previous theorem we will need the following definition and theorem which are due to Lauzon and Treil.
%%%%%%%%%%%%%%%%%%%%%%%%%%%%%%%%%%%%%%%%%%%%%%%%%%%%%%%%%%%%%%%%%%%%%

\begin{definition}\cite[p. 510]{Lauzon}
Let $X$ be a Hilbert space. 
\begin{itemize}
\item The upper linear codimension of a subset $K$ of $X$ is defined as 
$$\inf\left\{\mathrm{codim}(L)\,|\,L \mathrm{\; is\; a\; linear\; subspace\; of\; }K\right\}.$$ 
\item If $M$ and $N$ are subspaces of $X$ and $\varepsilon>0$ then we define the cone
$$\mathcal{K}^{\,\varepsilon}_M=\left\{x\in M\,|\,\mathrm{dist}(x,N)\leq \varepsilon\,\| x\|\right\}.$$
\end{itemize}
\end{definition}
%%%%%%%%%%%%%%%%%%%%%%%%%%%%%%%%%%%%%%%%%%%%%%%%%%%%%%%%%%%%%%%%%%%%%

\begin{theorem}\cite[Theorem 5.1]{Lauzon}
\label{geometric}
The subspaces $M$ and $N$ of a Hilbert space $X$ have a common complement if and only if for some small $\varepsilon>0$
the upper linear codimensions of the cones $\mathcal{K}^{\,\varepsilon}_M$ in $M$ and $\mathcal{K}^{\,\varepsilon}_N$ in $N$ coincide.
\end{theorem}
%%%%%%%%%%%%%%%%%%%%%%%%%%%%%%%%%%%%%%%%%%%%%%%%%%%%%%%%%%%%%%%%%%%%%

As Lauzon and Treil note, the equality of upper linear codimensions in the previous theorem can be thought of as an $\varepsilon$-analogue of the equality of dimensions (\ref{dimensions}). Now we can give the proof of Theorem \ref{asymptotic}.
%%%%%%%%%%%%%%%%%%%%%%%%%%%%%%%%%%%%%%%%%%%%%%%%%%%%%%%%%%%%%%%%%%%%%

\begin{proof}
Assume that $M$ and $N$ have a common complement and are not equivalently positioned. As we have already mentioned in the Introduction, condition (\ref{Lauzon}) then implies that
$$\mathrm{dim}(\mathcal{E}((0,1-\varepsilon))(M))=\infty$$
and this in turn implies that the operator $I-G^*G:M\rightarrow M$ is not compact (see \cite[Remark 0.5]{Lauzon}). But $$I-G^*G=\left. P_M(I-P_N)\right|_M$$ 
and hence $\left. (I-P_N)\right|_M$ is not compact. This implies that
$(I-P_N)(M)$ contains an infinite dimensional closed subspace and therefore by Remark \ref{asym} the subspace $M$ is not completely asymptotic to $N$. By observing that $(I-G^*G)|_M$ is not compact if and only if $(I-GG^*)|_N$ is not compact, see \cite[p. 502]{Lauzon}, we get that $N$ is also not completely asymptotic to $M$ and hence the required result follows.

Conversely, if $M$ and $N$ are equivalently positioned we are done by (\ref{Lauzon}). If this is not the case then $M$ and $N$ are not completely asymptotic to one another and thus there exists an infinite dimensional subspace $M_1$ of $M$ disjoint from $N$ with $M_1+N$ closed. Hence there exists $\varepsilon_1>0$ such that
\begin{equation}
\label{upperlinear}
M_1\cap \mathcal{K}^{\,\varepsilon_1}_M=\left\{0\right\}.
\end{equation}
It follows from (\ref{upperlinear}), using \cite[Lemma 5.12]{Schechter}, that $\mathcal{K}^{\,\varepsilon_1}_M$ does not contain a finite codimensional subspace of $M$ and thus the upper linear codimension of $\mathcal{K}^{\,\varepsilon_1}_M$ in $M$ is infinite. Similarly we can find $\varepsilon_2>0$ such that the upper linear codimension of $\mathcal{K}^{\,\varepsilon_2}_N$ in $N$ is infinite. Choosing $\varepsilon=\min\left\{\varepsilon_1,\,\varepsilon_2\right\}$ we get that the upper linear codimensions of $\mathcal{K}^{\,\varepsilon}_M$ in $M$ and of $\mathcal{K}^{\,\varepsilon}_N$ in $N$ coincide, since the space is separable, and hence, by Theorem \ref{geometric}, we get that $M$ and $N$ have a common complement.
\end{proof}
%%%%%%%%%%%%%%%%%%%%%%%%%%%%%%%%%%%%%%%%%%%%%%%%%%%%%%%%%%%%%%%%%%%%%

\begin{remark}
\label{corollary11aa}
An immediate corollary of the previous theorem is that in a separable Hilbert space, if two subspaces have closed sum, then the existence of a common complement is a matter of dimensions. More precisely let $X$ be a separable Hilbert space and $M$ and $N$ be subspaces of $X$ such that $M+N$ is closed. Then
$M$ and $N$ have a common complement in $X$ if and only if
$$\mathrm{dim} (M\ominus (M\cap N))=\mathrm{dim} (N\ominus (M\cap N)).$$ 
\end{remark}
%%%%%%%%%%%%%%%%%%%%%%%%%%%%%%%%%%%%%%%%%%%%%%%%%%%%%%%%%%%%%%%%%%%%%

We will now extend the result of the above remark to Banach spaces. To do that we will use Proposition \ref{proposition10}. It turns out again that for subspaces with closed sum everything works more or less as in the finite dimensional case.
%%%%%%%%%%%%%%%%%%%%%%%%%%%%%%%%%%%%%%%%%%%%%%%%%%%%%%%%%%%%%%%%%%%%%

\begin{proposition}
\label{proposition11}
Let $X$ be a Banach space and $M$ and $N$ be subspaces of $X$ such that $M+N$ is closed and $M\cap N=\{ 0\}$.
Then $M$ and $N$ have a common complement in $M\oplus N$ if and only if they are isomorphic.
\end{proposition}
%%%%%%%%%%%%%%%%%%%%%%%%%%%%%%%%%%%%%%%%%%%%%%%%%%%%%%%%%%%%%%%%%%%%

\begin{proof}
Suppose that $M$ and $N$ are isomorphic and let $U:M\rightarrow N$ be an isomorphism onto. Since $M+N$ is closed and $M\cap N=\{ 0\}$, there exists $c\geq 1$ such that
$$\| x\| +\| y\|\leq c\| x+y\|,$$
for all $x\in M$ and $y\in N$. For $C=c\;\mathrm{max}\left\{\| U\|,\,1\right\}\geq 1$ we have
$$
\| Ux+y\|\leq C\;\| x+y\|\,,
$$
for all $x\in M,\,y\in N$. 
Thus, by Proposition \ref{proposition10}, $M$ and $N$ have a common complement in $\overline{M+N}=M\oplus N$. The other direction follows immediately from Corollary \ref{corollary6}.
\end{proof}
%%%%%%%%%%%%%%%%%%%%%%%%%%%%%%%%%%%%%%%%%%%%%%%%%%%%%%%%%%%%%%%%%%%%

\begin{remark}
\label{remark23}
The result of Proposition \ref{proposition11} is not true if $M\cap N\ne\{ 0\}$. To see that let $X$ be a Banach space, $M$ be a subspace of $X$ and $N$ be a proper subspace of $M$ which is isomorphic to $M$. Then $M+N=M$ is closed, $M$ and $N$ are isomorphic, $M\cap N=N\ne\{ 0\}$ and $M$ and $N$ do not have a common complement. 
\end{remark}
%%%%%%%%%%%%%%%%%%%%%%%%%%%%%%%%%%%%%%%%%%%%%%%%%%%%%%%%%%%%%%%%%%%%

Combining Proposition \ref{proposition11} with Corollaries \ref{corollary1}, \ref{corollary3} and \ref{corollary10} we get the following:
%%%%%%%%%%%%%%%%%%%%%%%%%%%%%%%%%%%%%%%%%%%%%%%%%%%%%%%%%%%%%%%%%%%%

\begin{corollary}
\label{corollary11}
Let $X$ be a Banach space and $M$ and $N$ be subspaces of $X$ such that $M+N$ and  $M\cap N$ are complemented in $X$. 
If $P:X\rightarrow X$ is a bounded projection with $P(X)=M\cap N$, then $M$ and $N$ have a common complement in $X$ if and only if $(I-P)(M)$ and $(I-P)(N)$ are isomorphic. 
\end{corollary}
%%%%%%%%%%%%%%%%%%%%%%%%%%%%%%%%%%%%%%%%%%%%%%%%%%%%%%%%%%%%%%%%%%%%

In particular for Hilbert spaces we have the following generalization of Remark \ref{corollary11aa}. The separability of $X$ is no longer necessary. 
%%%%%%%%%%%%%%%%%%%%%%%%%%%%%%%%%%%%%%%%%%%%%%%%%%%%%%%%%%%%%%%%%%%%

\begin{corollary}
\label{corollary11a}
Let $X$ be a Hilbert space and $M$ and $N$ be subspaces of $X$ such that $M+N$ is closed. Then
$M$ and $N$ have a common complement in $X$ if and only if
$$\mathrm{dim} (M\ominus (M\cap N))=\mathrm{dim} (N\ominus (M\cap N)).$$ 
\end{corollary}
%%%%%%%%%%%%%%%%%%%%%%%%%%%%%%%%%%%%%%%%%%%%%%%%%%%%%%%%%%%%%%%%%%%%%

\begin{remark}
\label{remark11a}
In \cite[Proposition 2.2]{Lauzon} Lauzon and Treil proved that if $\left\|\left. P_N\right|_M\right\|<1$ and $\mathrm{dim} (M\ominus (M\cap N))=\mathrm{dim} (N\ominus (M\cap N))$, then $M$ and $N$ have a common complement. This follows immediately from the previous corollary, since $\left\|\left. P_N\right|_M\right\|<1$ implies that $M+N$ is closed (the converse is not in general true). 
\end{remark}
%%%%%%%%%%%%%%%%%%%%%%%%%%%%%%%%%%%%%%%%%%%%%%%%%%%%%%%%%%%%%%%%%%%%%

Even though the existence of a common complement does not imply that the subspaces are not completely asymptotic to one another, subspaces with a common complement share a similar property: If two subspaces $M$ and $N$ of a Hilbert space have a common complement then for any infinite dimensional subspace of $M$, with infinite codimension in $M$, we can find an infinite dimensional subspace of $N$, with infinite codimension in $N$, such that those two subspaces have a closed sum. 
%%%%%%%%%%%%%%%%%%%%%%%%%%%%%%%%%%%%%%%%%%%%%%%%%%%%%%%%%%%%%%%%%%%%%

\begin{proposition}
\label{proposition13}
Let $X$ be a Hilbert space and $M$ and $N$ be subspaces of $X$. If $M$ and $N$ have a common complement in $X$, then, for each subspace $M_1$ of $M$, with $\mathrm{dim}(M_1)=\mathrm{codim}_M(M_1)=\infty$, there exists a subspace $N_1$ of $N$, with $\mathrm{dim}(N_1)=\mathrm{codim}_N(N_1)=\infty$, such that $M_1+N_1$ is closed.
\end{proposition}
%%%%%%%%%%%%%%%%%%%%%%%%%%%%%%%%%%%%%%%%%%%%%%%%%%%%%%%%%%%%%%%%%%%%

\begin{proof}
Since $M$ and $N$ have a common complement in $X$, by Proposition \ref{proposition10}, there exist an isomorphism onto $U:N\rightarrow M$ and $C\geq 1$ such that 
\begin{equation}
\label{inequality1}
\| Ux+y\|\leq C\| x+y\| ,
\end{equation}
for all $x\in N$ and $y\in M$. Let $M_1$ be a subspace of $M$ with
$$\mathrm{dim}(M_1)=\mathrm{codim}_M(M_1)=\infty$$
and let $M_2=M_1^{\perp}$. Since $M_1\perp M_2$,
\begin{equation}
\label{inequality2}
\| y+z\|\geq \frac{1}{\sqrt{ 2}}\,(\| y\| +\| z\| ),
\end{equation}
for all $y\in M_1$ and $z\in M_2$. Let $N_1=U^{-1}(M_2)$. Then $N_1$ is a subspace of $N$ with
$$
\begin{array}{ccc}
\mathrm{dim}(N_1)=\mathrm{dim}(M_2)=\infty&
\mathrm{and}&
\mathrm{codim}_N(N_1)=\mathrm{dim}(M_1)=\infty .
\end{array}
$$
From inequalities (\ref{inequality1}), (\ref{inequality2}) and the fact that $U$ is an isomorphism we get that for $C'=\displaystyle\frac{1}{\sqrt{ 2}\,C}\;\mathrm{min}\left\{ \frac{1}{\| U^{-1}\|},1\right\}$, we have
$$
\| x+y\|\geq C'\;(\| x\| +\| y\| ),
$$
for all $x\in N_1$ and $y\in M_1$, and thus $M_1+N_1$ is closed.
\end{proof}
%%%%%%%%%%%%%%%%%%%%%%%%%%%%%%%%%%%%%%%%%%%%%%%%%%%%%%%%%%%%%%%%%%%%

\begin{remark}
The above proposition should not be misleading in the sense that the above construction  depends heavily on the fact that $\mathrm{codim}_M(M_1)=\infty$. If this is not the case then we cannot get a nontrivial conclusion (take for example $M_1=M$).
\end{remark}
%%%%%%%%%%%%%%%%%%%%%%%%%%%%%%%%%%%%%%%%%%%%%%%%%%%%%%%%%%%%%%%%%%%%

\section{A characterization via graphs}
\label{graphs}
%%%%%%%%%%%%%%%%%%%%%%%%%%%%%%%%%%%%%%%%%%%%%%%%%%%%%%%%%%%%%%%%%%%%

An idea which is widely used in the study of the relative position of a pair of subspaces is to represent it as a pair $\{\mathrm{Gr}(T),\mathrm{Gr}(S)\}$ of graphs of bounded or unbounded linear operators $T$ and $S$. The properties of those operators characterize the relative position of the pair. This idea goes back to Halmos \cite{Halmos}. The characterization of topologically complementary subspaces of Theorem \ref{Longstaff} is such a result. We shall use this approach to characterize pairs of subspaces with a common complement. 
%%%%%%%%%%%%%%%%%%%%%%%%%%%%%%%%%%%%%%%%%%%%%%%%%%%%%%%%%%%%%%%%%%%%

Recall the following characterizations of pairs of subspaces of a Hilbert space which are in generic position; the equivalence of (1) and (2) is \cite[Theorem 3]{Halmos} and of (1) and (3) is \cite[Theorem 1]{Halmos} (see also \cite[Lemma 4.1]{Araki}). If $\{ M, N\}$ is a pair of subspaces in a Hilbert space $X$, then the following are equivalent:
\begin{enumerate}
\item The pair $\{M, N\}$ is in generic position, i.e.
$$M\cap N=M\cap N^{\perp}=M^{\perp}\cap N=M^{\perp}\cap N^{\perp}=\{ 0\} .$$
\item There exist a Hilbert space $X_1$ and $T:X_1\rightarrow X_1$ positive injective contraction with $I-T$ injective, such that $\{ M, N\}$ is unitarily equivalent to $\{\mathrm{Gr}(-T), \mathrm{Gr}(T)\}$.  
\item There exist a Hilbert space $Y_1$ and an injective densely defined closed linear operator $S:D(S)\rightarrow Y_1$ with dense range, such that $\{ M, N\}$ is unitarily equivalent to $\{ Y_1\oplus\{ 0\}, \mathrm{Gr}(S)\}$.  
\end{enumerate}
We should note here that the equivalence of (1) and (2) is true even if $T$ is not a contraction. The equivalence of (1) and (2) has been extended by Longstaff and Panaia \cite[p.\ 3022]{Longstaff} (this result was previously stated without proof in \cite[p.\ 1158]{Papadakis}) as follows: If $\{ M, N\}$ is a pair of subspaces in a Hilbert space $X$, then the following are equivalent:
\begin{enumerate}
\item The pair $\{ M, N\}$ is in generalized generic position, i.e.
$$
\begin{array}{ccc}
M\cap N=M^{\perp}\cap N^{\perp}=\{ 0\}&
\mathrm{ and }&
\mathrm{dim}(M\cap N^{\perp})=\mathrm{dim}(M^{\perp}\cap N).
\end{array}
$$
\item There exist a Hilbert space $X_1$ and $T:X_1\rightarrow X_1$ positive injective contraction, such that the pair $\{ M, N\}$ is unitarily equivalent to $\{\mathrm{Gr}(-T), \mathrm{Gr}(T)\}$.
\end{enumerate}
Again we note that the equivalence of (1) and (2) is true even if $T$ is not a contraction.
%%%%%%%%%%%%%%%%%%%%%%%%%%%%%%%%%%%%%%%%%%%%%%%%%%%%%%%%%%%%%%%%%%%%

In the proposition that follows we extend the above results to obtain characterizations of subspaces in position $\mathrm{p}\,'$ (see the definition below) and of equivalently positioned subspaces via graphs, which we will use later.
%%%%%%%%%%%%%%%%%%%%%%%%%%%%%%%%%%%%%%%%%%%%%%%%%%%%%%%%%%%%%%%%%%%%

\begin{proposition}
\label{proposition14}
Let  $X$ be a Hilbert space and $\{ M, N\}$ be a pair of subspaces of $X$. Then the following hold:
\begin{enumerate}
\item The following are equivalent:
\begin{enumerate}
\item The pair $\{ M, N\}$ is in position $\mathit{p}\,'$, i.e.
$$M\cap N^{\perp}=M^{\perp}\cap N=\{ 0\} .$$
\item There exist Hilbert spaces $X_1$, $X_2$ and a contraction $$T:X_1\rightarrow X_2$$ with 
$I-T^*T$ injective, such that $\{ M, N\}$ is unitarily equivalent to 
$\{\mathrm{Gr}(-T), \mathrm{Gr}(T)\}$.
\item There exist Hilbert spaces $Y_1$ and $Y_2$ and a densely defined closed linear operator $S:D(S)\rightarrow Y_2$, such that $\{ M, N\}$ is unitarily equivalent to $\{ Y_1\oplus\{ 0\}, \mathrm{Gr}(S)\}$.
\end{enumerate}
\item The following are equivalent:
\begin{enumerate}
\item The pair $\{ M, N\}$ is equivalently positioned.
\item There exist Hilbert spaces $X_1$ and $X_2$ and a contraction $$T:X_1\rightarrow X_2,$$ such that $\{ M, N\}$ is unitarily equivalent to $\{\mathrm{Gr}(-T), \mathrm{Gr}(T)\}$.
\end{enumerate}
\end{enumerate}
\end{proposition}
%%%%%%%%%%%%%%%%%%%%%%%%%%%%%%%%%%%%%%%%%%%%%%%%%%%%%%%%%%%%%%%%%%%%

\begin{proof}
(1) (a)$\Rightarrow$(b): Let 
$$
\begin{array}{ccc}
M_1=M\ominus (M\cap N)&
\mathrm{and}&
N_1=N\ominus(M\cap N).
\end{array}
$$ 
It is easy to see that, since $\{ M, N\}$ is in position $\mathrm{p}\,'$, $\{ M_1, N_1\}$ is in generic position in $\overline{M_1+N_1}$. Thus, by the first of Halmos' results that we mentioned above, there exist a Hilbert space $Z$ and $A:Z\rightarrow Z$ positive injective contraction with $I-A$ injective, such that the pair $\{ M_1, N_1\}$ is unitarily equivalent to the pair $\{\mathrm{Gr}(-A), \mathrm{Gr}(A)\}$, via an isometry $V$ from $\overline{M_1+N_1}$ onto $Z\oplus Z$. 

Let 
$$X_1=Z\oplus (M\cap N),\;\;X_2=Z\oplus (\overline{M+N})^{\perp},$$
$T:X_1\rightarrow X_2$, with $T((z,y))=(Az,0)$, for all $z\in Z$ and $y\in M\cap N$, and 
$$U:X\rightarrow X_1\oplus X_2,$$ 
with 
$$Ux=((P_1VP_{\overline{M_1+N_1}}\,x,P_{M\cap N}x), (P_2VP_{\overline{M_1+N_1}}\,x,P_{(\overline{M+N})^{\perp}}x)),$$ 
for all $x\in X$, where $P_1$ and $P_2$ are the projections from $Z\oplus Z$ onto the first and second component respectively. It is easy to see that $T$ is a contraction with $I-T^*T$ injective and that $U$ is an isometry onto $X_1\oplus X_2$ taking $\{ M, N\}$ to $\{\mathrm{Gr}(-T), \mathrm{Gr}(T)\}$.  

\noindent (b)$\Rightarrow$(a): A straightforward calculation shows that
\begin{equation}
\label{equation1}
\begin{array}{c}
\mathrm{Gr}(-T)\cap\mathrm{Gr}(T)^{\perp}=\mathrm{Gr}\left(-T|_{\mathrm{Ker}(I-T^*T)}\right)\\[5pt]
\mathrm{Gr}(-T)^{\perp}\cap\mathrm{Gr}(T)=\mathrm{Gr}\left(T|_{\mathrm{Ker}(I-T^*T)}\right).
\end{array}
\end{equation}
Since $\mathrm{Ker}(I-T^*T)=\{ 0\}$,
$$\mathrm{Gr}(-T)\cap\mathrm{Gr}(T)^{\perp}=\{ 0\}=\mathrm{Gr}(-T)^{\perp}\cap\mathrm{Gr}(T)$$
and so $\{\mathrm{Gr}(-T), \mathrm{Gr}(T)\}$ is in position $\mathrm{p}\,'$. Since position $\mathrm{p}\,'$ is preserved under unitary equivalence, the pair $\{ M, N\}$ is also in position $\mathrm{p}\,'$.

\noindent (a)$\Rightarrow$(c): It follows from the second of Halmos' results that we mentioned above in exactly the same manner as (a)$\Rightarrow$(b).

\noindent (c)$\Rightarrow$(a): Since $S$ is densely defined and closed, $S^*:D(S^*)\rightarrow Y_1$ is also densely defined and closed. Thus
$$(Y_1\oplus\{0\})\cap\mathrm{Gr}(S)^{\perp}=(Y_1\oplus\{0\})\cap\mathrm{Gr}(-S^*)=\{ 0\}$$
and
$$(Y_1\oplus\{0\})^{\perp}\cap\mathrm{Gr}(S)=(\{0\}\oplus Y_2)\cap\mathrm{Gr}(S)=\{ 0\}.$$

\noindent (2) (a)$\Rightarrow$(b): It follows from the characterization of subspaces in generalized generic position due to Longstaff and Panaia that we mentioned above in exactly the same manner as (a)$\Rightarrow$(b) in (1).

\noindent (b)$\Rightarrow$(a): By Equation (\ref{equation1}),
$$\mathrm{dim}(\mathrm{Gr}(-T)\cap\mathrm{Gr}(T)^{\perp})=\mathrm{dim}(\mathrm{Ker}(I-T^*T))=
\mathrm{dim}(\mathrm{Gr}(-T)^{\perp}\cap\mathrm{Gr}(T))$$
and thus the pair $\{\mathrm{Gr}(-T), \mathrm{Gr}(T)\}$ is equivalently positioned. Therefore, since the pair $\{ M, N\}$ is unitarily equivalent to $\{\mathrm{Gr}(-T), \mathrm{Gr}(T)\}$, it is also equivalently positioned.
\end{proof}
%%%%%%%%%%%%%%%%%%%%%%%%%%%%%%%%%%%%%%%%%%%%%%%%%%%%%%%%%%%%%%%%%%%%

As before we have that the equivalence of (a) and (b) both in (1) and in (2) holds even if $T$ is not a contraction.
%%%%%%%%%%%%%%%%%%%%%%%%%%%%%%%%%%%%%%%%%%%%%%%%%%%%%%%%%%%%%%%%%%%%

We move now to the characterization via graphs of pairs of subspaces of a Banach space which have a common complement. As we already said two subspaces of a Banach space have a common complement if and only if there exists a bounded projection onto one of them, the restriction of which on to the other is an isomorphism onto. The following form of that result, which specifies the projection, can be found in \cite[Lemma 2.1(a)]{Dirr}. 
%%%%%%%%%%%%%%%%%%%%%%%%%%%%%%%%%%%%%%%%%%%%%%%%%%%%%%%%%%%%%%%%%%%%

\begin{lemma}
\label{proposition7}
Let $X$ be a Banach space and $M$, $N$ and $K$ be subspaces of $X$. The following are equivalent:
\begin{enumerate}
\item $K$ is a common complement of $M$ and $N$ in $X$.
\item $X=M\oplus K$ and $\left. P_{M\| K}\right|_N:N\rightarrow M$ is an isomorphism onto.
\end{enumerate}
Moreover if things are as above, then
$$
\left( \left. P_{M\| K}\right|_N\right)^{-1}=\left. P_{N\| K}\right|_M.
$$
\end{lemma}
%%%%%%%%%%%%%%%%%%%%%%%%%%%%%%%%%%%%%%%%%%%%%%%%%%%%%%%%%%%%%%%%%%%%

Our main result for this section is the following:
%%%%%%%%%%%%%%%%%%%%%%%%%%%%%%%%%%%%%%%%%%%%%%%%%%%%%%%%%%%%%%%%%%%%

\begin{theorem}
\label{theorem1}
Let $X$ be a Banach space and $M$ and $N$ be subspaces of $X$. Then the following are equivalent:
\begin{enumerate}
\item $M$ and $N$ have a common complement in $X$.
\item There exist Banach spaces $X_1$, $X_2$ and $T,\,S:X_1\rightarrow X_2$ bounded linear operators such that $\{ M, N\}$ is isomorphic to $\{\mathrm{Gr}(T),\mathrm{Gr}(S)\}$.
\end{enumerate}
\end{theorem}
%%%%%%%%%%%%%%%%%%%%%%%%%%%%%%%%%%%%%%%%%%%%%%%%%%%%%%%%%%%%%%%%%%%%

\begin{proof}
(1)$\Rightarrow$(2): Let $K$ be a common complement of $M$ and $N$ in $X$ and let $K'$ be a complement of $K$ in $X$. Then
$$X=M\oplus K=N\oplus K=K'\oplus K$$
and thus $K$ is a common complement of $M$ and $K'$ and of $N$ and $K'$. If $P=P_{K'\| K}$ then, by Lemma \ref{proposition7}, the operators
$$
\begin{array}{ccc}
G_1=P|_M:M\rightarrow K'&
\mathrm{and}&
G_2=P|_N:N\rightarrow K'
\end{array}
$$
are isomorphisms onto. Let $X_1=K'$, $X_2=K$ and define 
$$T:X_1\rightarrow X_2$$
with $Tx=G_1^{-1}x-x$, for all $x\in K'$, and 
$$S:X_1\rightarrow X_2$$
with $Sx=G_2^{-1}x-x$, for all $x\in K'$. Also let
$$U:X\rightarrow X_1\oplus X_2$$
with $Ux=(Px,(I-P)x)$, for all $x\in X$. By Lemma \ref{proposition7}, $G_1^{-1}=\left. P_{M\| K}\right|_{K'}$. Thus, for all $x\in K'$,
$$Tx=G_1^{-1}x-x=\left. P_{M\| K}\right|_{K'}x-x=-\left. P_{K\| M}\right|_{K'}x$$
and therefore $T$ is a well-defined bounded linear operator. Similarly we get that $S$ is a well-defined bounded linear operator. Since $X=K'\oplus K$, $U$ is an isomorphism onto. We will show that $U(M)=\mathrm{Gr}(T)$. To this end let $x\in M$. Then
\begin{eqnarray*}
Ux&=&(Px,(I-P)x)\\
&=&(G_1x,x-G_1x)\\
&=&(G_1x,G_1^{-1}G_1x-G_1x)\\
&=&(G_1x,TG_1x)\in \mathrm{Gr}(T).
\end{eqnarray*}
On the other hand if $(x,Tx)\in\mathrm{Gr}(T)$, then, since $G_1$ is an isomorphism onto, there exists
$y\in M$ such that $G_1y=x$ and hence
\begin{eqnarray*}
(x,Tx)&=&(G_1y,TG_1y)\\
&=&(G_1y,G_1^{-1}G_1y-G_1y)\\
&=&(G_1y,y-G_1y)\\
&=&(Py,(I-P)y)\in U(M).
\end{eqnarray*}
Thus
$$U(M)=\mathrm{Gr}(T).$$
Similarly we can prove that
$$U(N)=\mathrm{Gr}(S).$$
Therefore the pair $\{ M,N\}$ is isomorphic to the pair $\{\mathrm{Gr}(T),\mathrm{Gr}(S)\}$.

\noindent (2)$\Rightarrow$(1): It is easy to see that $\left\{0\right\}\oplus X_2$ is a common complement of 
$\mathrm{Gr}(T)$ and $\mathrm{Gr}(S)$ in $X_1\oplus X_2$. Therefore, by Proposition \ref{proposition5}, $M$ and $N$ have a common complement in $X$.
\end{proof}
%%%%%%%%%%%%%%%%%%%%%%%%%%%%%%%%%%%%%%%%%%%%%%%%%%%%%%%%%%%%%%%%%%%%

\begin{remark}
\label{remark26}
It is obvious that the equivalence of (1) and (2) still holds if we replace the pair of subspaces with a family of subspaces.
\end{remark}
%%%%%%%%%%%%%%%%%%%%%%%%%%%%%%%%%%%%%%%%%%%%%%%%%%%%%%%%%%%%%%%%%%%%

If $X$ is a Hilbert space, then in the proof of (1)$\Rightarrow$(2) of Theorem \ref{theorem1} we can take $K'=K^{\perp}$. In that case $U$ is an isometry and so we get the following improvement of Theorem \ref{Longstaff}.
%%%%%%%%%%%%%%%%%%%%%%%%%%%%%%%%%%%%%%%%%%%%%%%%%%%%%%%%%%%%%%%%%%%%

\begin{corollary}
\label{theorem3}
Let $X$ be a Hilbert space and $M$ and $N$ be subspaces of $X$. Then the following are equivalent:
\begin{enumerate}
\item $M$ and $N$ have a common complement in $X$.
\item There exist Hilbert space $X_1$, $X_2$ and $T,\,S:X_1\rightarrow X_2$ bounded linear operators such that $\{ M, N\}$ is unitarily equivalent to $\{\mathrm{Gr}(T),\mathrm{Gr}(S)\}$.
\end{enumerate}
\end{corollary}
%%%%%%%%%%%%%%%%%%%%%%%%%%%%%%%%%%%%%%%%%%%%%%%%%%%%%%%%%%%%%%%%%%%%

\begin{remark}
\label{theorem3a}
If we replace unitarily equivalent with similar in (2) of the previous corollary, then we can prove the result combining the equivalence of (1) and (2) in Theorem \ref{Longstaff} and the comments following Corollary \ref{corollary3}. In a similar manner using the equivalence of (1) and (3) in Theorem \ref{Longstaff} we get the following characterization. Let $X$ be a Hilbert space, $M$ and $N$ be subspaces of $X$ and $M_1=M\ominus (M\cap N)$, $N_1=N\ominus (M\cap N)$. Then the following are equivalent:
\begin{enumerate}
\item $M$ and $N$ have a common complement in $X$.
\item The pair $\{ M_1, N_1\}$ in $\overline{M_1+N_1}$ is similar to a pair of subspaces in generic position.
\end{enumerate}
\end{remark}
%%%%%%%%%%%%%%%%%%%%%%%%%%%%%%%%%%%%%%%%%%%%%%%%%%%%%%%%%%%%%%%%%%%%

If in the proof of (1)$\Rightarrow$(2) of Theorem \ref{theorem1} we take $K'=M$ we get that $\{ M,N\}$ is isomorphic to $\{ X_1\oplus\{ 0\},\mathrm{Gr}(S)\}$. Moreover the pair $\{ X_1\oplus\{ 0\},\mathrm{Gr}(S)\}$ is isomorphic to the pair $\left\{ \mathrm{Gr}\left(-\frac{S}{2}\right),\mathrm{Gr}\left(\frac{S}{2}\right)\right\}$ via the isomorphism
$$
\left[
\begin{array}{cc}
I_{X_1}&0\\
-\frac{S}{2}&I_{X_2}
\end{array}
\right]
:X_1\oplus X_2\rightarrow X_1\oplus X_2.
$$
Thus we get two more conditions equivalent to the existence of a common complement.
%%%%%%%%%%%%%%%%%%%%%%%%%%%%%%%%%%%%%%%%%%%%%%%%%%%%%%%%%%%%%%%%%%%%

\begin{corollary}
\label{theorem1a}
Let $X$ be a Banach space and $M$ and $N$ be subspaces of $X$. Then the following are equivalent:
\begin{enumerate}
\item $M$ and $N$ have a common complement in $X$.
\item There exist Banach spaces $X_1$, $X_2$ and $S:X_1\rightarrow X_2$ bounded linear operator such that $\{ M, N\}$ is isomorphic to $\{ X_1\oplus\{ 0\},\mathrm{Gr}(S)\}$.
\item There exist Banach spaces $Y_1$, $Y_2$ and $T:Y_1\rightarrow Y_2$ bounded linear operator such that $\{ M, N\}$ is isomorphic to $\{ \mathrm{Gr}(-T),\mathrm{Gr}(T)\}$.
\end{enumerate}
\end{corollary}
%%%%%%%%%%%%%%%%%%%%%%%%%%%%%%%%%%%%%%%%%%%%%%%%%%%%%%%%%%%%%%%%%%%%

It is easy to see that a pair $\{ M, N\}$ of subspaces of a Hilbert space $X$ is in position $\mathrm{p}\,'$ if and only if $M^{\perp}$ is a common topological complement of $M$ and $N$ in $X$ (or equivalently $N^{\perp}$ is a common topological complement of $M$ and $N$ in $X$). In the following corollary we characterize pairs $\{ M, N\}$ of subspaces of a Hilbert space $X$ for which $M^{\perp}$ is a common complement (or equivalently $N^{\perp}$ is a common complement). This characterization can be thought of as an analogue of Proposition \ref{proposition14}(1).
%%%%%%%%%%%%%%%%%%%%%%%%%%%%%%%%%%%%%%%%%%%%%%%%%%%%%%%%%%%%%%%%%%%%

\begin{corollary}
\label{theorem1b}
Let $X$ be a Hilbert space, $M$ and $N$ be subspaces of $X$ and $M_1=M\ominus (M\cap N)$, $N_1=N\ominus (M\cap N)$. Then the following are equivalent:
\begin{enumerate}
\item $M^{\perp}$ is a common complement of $M$ and $N$ in $X$.
\item There exist Hilbert spaces $X_1$, $X_2$ and an injective bounded linear operator $T:X_1\rightarrow X_2$ with $I-T^*T$ injective and $I-TT^*$ onto such that $\{ M, N\}$ is unitarily equivalent to $\{\mathrm{Gr}(-T),\mathrm{Gr}(T)\}$.
\item There exist Hilbert spaces $Y_1$, $Y_2$ and $S:Y_1\rightarrow Y_2$ bounded linear operator such that $\{ M, N\}$ is unitarily equivalent to $\{ Y_1\oplus\{ 0\},\mathrm{Gr}(S)\}$.
\end{enumerate}
\end{corollary}
%%%%%%%%%%%%%%%%%%%%%%%%%%%%%%%%%%%%%%%%%%%%%%%%%%%%%%%%%%%%%%%%%%%%

\begin{proof}
(1)$\Rightarrow$(2): From (a)$\Rightarrow$(b) of Proposition \ref{proposition14}(1), there exist Hilbert spaces $X_1$, $X_2$ and an injective bounded linear operator $T:X_1\rightarrow X_2$ with $I-T^*T$ injective, such that the pair $\{ M, N\}$ is unitarily equivalent to $\{\mathrm{Gr}(-T),\mathrm{Gr}(T)\}$. Moreover $\mathrm{Gr}(-T)^{\perp}=\mathrm{Gr}(T^*)$ must be a common complement of $\mathrm{Gr}(-T)$ and $\mathrm{Gr}(T)$. In particular we must have 
$$\mathrm{Gr}(T)+\mathrm{Gr}(T^*)=X_1\oplus X_2.$$
A straightforward calculation shows that 
\begin{equation}
\label{equation2}
\mathrm{Gr}(T)+\mathrm{Gr}(T^*)=\mathrm{Gr}(T)\oplus (\{ 0\}\oplus R(I-TT^*))
\end{equation}
and hence $I-TT^*$ is onto.

\noindent (2)$\Rightarrow$(1): As in the proof of (b)$\Rightarrow$(a) of Proposition \ref{proposition14}.(1) using (\ref{equation2}).

\noindent (1)$\Rightarrow$(3): Take $K=M^{\perp}$ and $K'=M$ in the proof of (1)$\Rightarrow$(2) in Theorem \ref{theorem1}.

\noindent (3)$\Rightarrow$(1): It is straightforward since
$$(Y_1\oplus\{ 0\})^{\perp}\oplus\mathrm{Gr}(S)=(\{ 0\}\oplus Y_2)\oplus\mathrm{Gr}(S)=Y_1\oplus Y_2.$$
\end{proof}
%%%%%%%%%%%%%%%%%%%%%%%%%%%%%%%%%%%%%%%%%%%%%%%%%%%%%%%%%%%%%%%%%%%%

We already mentioned in the Introduction that an immediate consequence of Theorem \ref{Treil} is that a pair of equivalently positioned subspaces always has a common complement (note that this can also be proved combining Theorem \ref{theorem1} and Proposition \ref{proposition14}(2)). Moreover in finite dimensions all pairs with a common complement are of this kind. The following example shows that this is not the case in infinite dimensions.
%%%%%%%%%%%%%%%%%%%%%%%%%%%%%%%%%%%%%%%%%%%%%%%%%%%%%%%%%%%%%%%%%%%%

\begin{example}
\label{example}
Let $X=l^2(\mathbb{N})$, $\left\{e_n\,|\,n\in\mathbb{N}\right\}$ be the standard orthonormal basis of $X$,
$$
\begin{array}{lcl}
f_n=\displaystyle{\frac{1}{2}}\,e_{2n-1}-\displaystyle{\frac{\sqrt{3}}{2}}\,e_{2n},\;\;n\geq 1,&&\\[10pt]
g_n=\displaystyle{\frac{1}{2}}\,e_{2n-1}+\displaystyle{\frac{\sqrt{3}}{2}}\,e_{2n},\;\;n\geq 1,&\mathrm{and}& g_0=e_0,
\end{array}
$$ 
$M=\overline{\mathrm{span}}\left\{f_n\,|\,n\geq 1\right\}$ and $N=\overline{\mathrm{span}}\left\{g_n\,|\,n\in\mathbb{N}\right\}$. 
Then
$$\mathrm{dim}(M\cap N^\perp)=0\neq 1=\mathrm{dim}(M^\perp\cap N).$$
On the other hand $\mathrm{dim}\,M=\mathrm{dim}\,N$, $M\cap N=\{0\}$ and $M+N$ is closed. Hence, by Corollary \ref{corollary11a}, $M$ and $N$ have a common complement in $X$.
\end{example}
%%%%%%%%%%%%%%%%%%%%%%%%%%%%%%%%%%%%%%%%%%%%%%%%%%%%%%%%%%%%%%%%%%%%

Nevertheless, up to isomorphism, all pairs with a common complement are equivalently positioned:
%%%%%%%%%%%%%%%%%%%%%%%%%%%%%%%%%%%%%%%%%%%%%%%%%%%%%%%%%%%%%%%%%%%%

\begin{proposition}
\label{theorem4}
Let $X$ be a Hilbert space and $M$ and $N$ be subspaces of $X$. Then the following are equivalent:
\begin{enumerate}
\item $M$ and $N$ have a common complement in $X$.
\item There exists a Hilbert space $X_1$ and a pair $\{ M_1, N_1\}$ of equivalently positioned subspaces of $X_1$ such that the pair $\{ M, N\}$ is similar to the pair $\{ M_1, N_1\}$.
\end{enumerate}
\end{proposition}
%%%%%%%%%%%%%%%%%%%%%%%%%%%%%%%%%%%%%%%%%%%%%%%%%%%%%%%%%%%%%%%%%%%%

\begin{proof}
It follows from (1)$\Leftrightarrow$(3) of Corollary \ref{theorem1a} and Proposition \ref{proposition14}.(2).
\end{proof}
%%%%%%%%%%%%%%%%%%%%%%%%%%%%%%%%%%%%%%%%%%%%%%%%%%%%%%%%%%%%%%%%%%%%

\begin{remark}
\label{theorem4a}
We can also prove the previous proposition using the characterization at the end of Remark \ref{theorem3a}.
\end{remark}
%%%%%%%%%%%%%%%%%%%%%%%%%%%%%%%%%%%%%%%%%%%%%%%%%%%%%%%%%%%%%%%%%%%%

\section{A characterization via involutions}
\label{involutions}
%%%%%%%%%%%%%%%%%%%%%%%%%%%%%%%%%%%%%%%%%%%%%%%%%%%%%%%%%%%%%%%%%%%%

Another way to study the relative position of a pair $\{ M, N\}$ of subspaces is to find an operator with ``nice geometric properties'' which exchanges $M$ and $N$. This idea can be traced back to Dixmier's work in \cite{Dixmier}. We recall that if $X$ is a Hilbert space, then we say that a linear operator $S:X\rightarrow X$ is a symmetry if $S$ is unitary and $S^2=I$. We remind the reader some characterizations of relative positions via symmetries which can either be found in \cite[Section II1]{Dixmier} and in \cite[Theorem 4.1]{Davis} or are straightforward corollaries of the results contained there (for generalizations of some of those results to an algebraic context see \cite[Corollary 1]{Holland}, \cite[Theorem 3]{Berberian1} and \cite[Theorem 1.1]{Maeda}).
%%%%%%%%%%%%%%%%%%%%%%%%%%%%%%%%%%%%%%%%%%%%%%%%%%%%%%%%%%%%%%%%%%%%

Let $X$ be a Hilbert space and $M$ and $N$ be subspaces of $X$. Then the following hold:
\begin{enumerate}
\item If the pair $\{ M, N\}$ is in generic position, then there exists a unique symmetry $S:X\rightarrow X$ with $\langle Sx,x\rangle >0$, for all $x\in M\setminus\{ 0\}$, such that $S(M)=N$.   
\item The pair $\{ M, N\}$ is in position $\mathrm{p}\,'$ and $M^{\perp}\cap N^{\perp}=\{ 0\}$ if and only if there exists a unique symmetry $S:X\rightarrow X$ with $\langle Sx,x\rangle >0$, for all $x\in M\setminus\{ 0\}$, such that $S(M)=N$.
\item The pair $\{ M, N\}$ is in position $\mathrm{p}\,'$ if and only if there exists a symmetry $S:X\rightarrow X$ with $\langle Sx,x\rangle >0$, for all $x\in M\setminus\{ 0\}$, such that $S(M)=N$.
\item The pair $\{ M, N\}$ is equivalently positioned if and only if there exists a symmetry $S:X\rightarrow X$ such that $S(M)=N$. Moreover $S$ can be chosen so that $\langle Sx,x\rangle \geq 0$, for all $x\in M$.
\end{enumerate}
Note that $S(M)=N$ immediately implies that $S(N)=M$ and that all the positivity conditions also hold for the elements of $N$. 
%%%%%%%%%%%%%%%%%%%%%%%%%%%%%%%%%%%%%%%%%%%%%%%%%%%%%%%%%%%%%%%%%%%%

Since for a Hilbert space the existence of a symmetry exchanging two subspaces is equivalent to them being equivalently positioned and subspaces with a common complement need not be equivalently positioned, symmetries are not well-suited for our problem. Moreover, since we want to work in Banach spaces, we must ``translate'' the positivity conditions. It is easy to see that $\langle Sx,x\rangle \geq 0$, for all $x\in M$, if and only if $\| x+Sx\|\geq\sqrt{2}\;\| x\|$, for all $x\in M$. We start by recalling the definition of an involution, since involutions will replace symmetries in our characterization.
%%%%%%%%%%%%%%%%%%%%%%%%%%%%%%%%%%%%%%%%%%%%%%%%%%%%%%%%%%%%%%%%%%%%

\begin{definition}
\label{definition2}
Let $X$ be a Banach space and $S:X\rightarrow X$ be a bounded linear operator. We will say that $S$ is an involution if $S^2=I$. If $S$ is also an isometry, then we will say that $S$ is a symmetry.
\end{definition}
%%%%%%%%%%%%%%%%%%%%%%%%%%%%%%%%%%%%%%%%%%%%%%%%%%%%%%%%%%%%%%%%%%%%

We can now prove our result.
%%%%%%%%%%%%%%%%%%%%%%%%%%%%%%%%%%%%%%%%%%%%%%%%%%%%%%%%%%%%%%%%%%%%

\begin{theorem}
\label{theorem2}
Let $X$ be a Banach space and $M$ and $N$ be subspaces of $X$ with $\overline{M+N}=X$. Then the following are equivalent:
\begin{enumerate}
\item $M$ and $N$ have a common complement in $X$.
\item There exist an involution $S:X\rightarrow X$ and $C>0$ with
$$\| x+Sx\|\geq C\,\| x\|,$$
for all $x\in M$, such that $S(M)=N$.
\end{enumerate}
\end{theorem}
%%%%%%%%%%%%%%%%%%%%%%%%%%%%%%%%%%%%%%%%%%%%%%%%%%%%%%%%%%%%%%%%%%%%

\begin{proof}
(1)$\Rightarrow$(2): Since $M$ and $N$ have a common complement in $X$, by Corollary \ref{theorem1a}, there exist Banach spaces $Y_1$, $Y_2$ and a bounded linear operator $T:Y_1\rightarrow Y_2$, such that $\{ M, N\}$ is isomorphic, via an isomorphism onto $U:X\rightarrow Y_1\oplus Y_2$, to $\{\mathrm{Gr}(-T),\mathrm{Gr}(T)\}$. In the rest of the proof consider $Y_1\oplus Y_2$ equipped with the $1$-norm. Let $S_1:Y_1\oplus Y_2\rightarrow Y_1\oplus Y_2$, with
$$S_1((x_1,x_2))=(x_1,-x_2),$$
for all $(x_1,x_2)\in Y_1\oplus Y_2$. Obviously $S_1$ is a symmetry that exchanges $\mathrm{Gr}(-T)$ and $\mathrm{Gr}(T)$. Moreover
$$
\| (x_1,-Tx_1)+S_1((x_1,-Tx_1))\|=\| (2x_1,0)\|\geq \frac{2}{1+\| T\|}\;\| (x_1,-Tx_1)\|,
$$
for all $(x_1,-Tx_2)\in\mathrm{Gr}(-T)$. Hence if $S=U^{-1}S_1U$ and
$$C=\frac{2}{\|U\|\|U^{-1}\|(1+\| T\|)},$$
then $S$ is an involution such that $S(M)=N$, with
$$
\| x+Sx\|\geq C\;\| x\|,
$$
for all $x\in M$.

\noindent (2)$\Rightarrow$(1): Since $S$ is an involution, if
$$
\begin{array}{ccc}
K'=\left\{ x\in X|\,Sx=x\right\}&
\mathrm{and}&
K=\left\{ x\in X|\,Sx=-x\right\},
\end{array}
$$
then it is well-known that
$$X=K\oplus K'$$
and  
$$P=P_{K'\| K}=\frac{1}{2}(I+S).$$
We shall prove that $K$ is a common complement of $M$ and $K'$. By Lemma \ref{proposition7}, it is enough to show that the operator
$$P|_M:M\rightarrow K'$$
is an isomorphism onto $K'$. To this end let $x\in M$. Then
$$
\| P|_Mx\|=\displaystyle\frac{1}{2}\;\| x+Sx\|\geq
\frac{C}{2}\;\| x\|
$$
and thus $P|_M$ is injective with closed range. On the other hand, by hypothesis, $M+N$ is dense in $X$ and thus $P(M+N)$ is dense in $K'$. Hence the set 
$$\left\{ x+Sx\;|\; x\in M+N\right\}$$ 
is also dense in $K'$. But if $x\in M+N$ then, since $S(M)=N$, there exist $z,\,w\in M$ such that $x=z+Sw$. Then
$$x+Sx=z+Sw+S(z+Sw)=(z+w)+S(z+w)$$ 
which implies that
$$\left\{ x+Sx\;|\; x\in M+N\right\}=\left\{ y+Sy\;|\; y\in M\right\}.$$
Therefore the set 
$$\left\{ y+Sy|\; y\in M\right\}$$
is also dense in $K'$. Since that set coincides with the range of $P|_M$, the range of $P|_M$ is a dense and closed subspace of $K'$ and so $P|_M$ is an isomorphism onto. Similarly, since $S(M)=N$ and $\| x+Sx\|\geq C\,\| x\|$, for all $x\in M$, immediately imply that $S(N)=M$ and $\| x+Sx\|\geq C\,\| x\|$, for all $x\in N$, we can show that $K$ is a common complement of $N$ and $K'$. Therefore $M$ and $N$ have a common complement in $X$.
\end{proof}
%%%%%%%%%%%%%%%%%%%%%%%%%%%%%%%%%%%%%%%%%%%%%%%%%%%%%%%%%%%%%%%%%%%%

\begin{remark}
\label{theorem2a}
Recall that a bounded linear operator $S:X\rightarrow X$ on a Banach space $X$ is called accretive if for every $x\in X$, there exists $x^\ast\in X^\ast$ such that $\langle x^\ast,x\rangle=\| x^\ast\|^2=\|x\|^2$ and $\langle x^\ast,Sx\rangle\geq 0$, where $\langle \cdot,\cdot\rangle$ is the duality product between $X^*$ and $X$. It is easy to see that if $S$ is accretive, then $\| x+Sx\|\geq \| x\|$, for all $x\in X$. Thus a particular class of pairs of subspaces of a Banach space which have a common complement consists of those pairs interchanged by an involution accretive on their union. It would be interesting to characterize this class and its subclass for which $S$ is a symmetry, the members of which are a Banach space version of equivalently positioned subspaces. 
\end{remark}
%%%%%%%%%%%%%%%%%%%%%%%%%%%%%%%%%%%%%%%%%%%%%%%%%%%%%%%%%%%%%%%%%%%%

We finish with a characterization of pairs $\{ M, N\}$ of subspaces of a Hilbert space for which $M^{\perp}$ is a common complement.
%%%%%%%%%%%%%%%%%%%%%%%%%%%%%%%%%%%%%%%%%%%%%%%%%%%%%%%%%%%%%%%%%%%%

\begin{proposition}
\label{theorem2b}
Let $X$ be a Hilbert space and $M$ and $N$ be subspaces of $X$. Then the following are equivalent:
\begin{enumerate}
\item $M^{\perp}$ is a common complement of $M$ and $N$ in $X$.
\item There exists a symmetry $S:X\rightarrow X$ and $C<1$ with $\langle Sx,x\rangle >0$, for all $x\in M\setminus\{ 0\}$, and $|\langle Sx,y\rangle|<C\|x\|\|y\|$, for all $x\in M$ and $y\in M^{\perp}$, such that $S(M)=N$.
\end{enumerate}
\end{proposition}
%%%%%%%%%%%%%%%%%%%%%%%%%%%%%%%%%%%%%%%%%%%%%%%%%%%%%%%%%%%%%%%%%%%%
\begin{proof}
By (4) in the beginning of this section, $M^{\perp}\cap N=\{ 0\}$ and $M^{\perp}+N$ is dense if and only if there exists a symmetry $S:X\rightarrow X$ with $\langle Sx,x\rangle >0$, for all $x\in M\setminus\{ 0\}$, such that $S(M)=N$. The result follows immediately since $M^{\perp}+N=M^{\perp}+S(M)$ is closed if and only if there exists $C<1$ such that
$$|\langle Sx,y\rangle|<C\|Sx\|\|y\|=C\|x\|\|y\|,$$
for all $x\in M$ and $y\in M^{\perp}$. 
\end{proof}
%%%%%%%%%%%%%%%%%%%%%%%%%%%%%%%%%%%%%%%%%%%%%%%%%%%%%%%%%%%%%%%%%%%%

\begin{Acknowledgments}
We would like to thank Prof. M. Anoussis and Prof. A. Katavolos for many helpful discussions.
\end{Acknowledgments}
%%%%%%%%%%%%%%%%%%%%%%%%%%%%%%%%%%%%%%%%%%%%%%%%%%%%%%%%%%%%%%%%%%%%

%%%%%%%%%%%%%%%%%%%%%%%%%%%%%%%%%%%%%%%%%%%%%%%%%%%%%%%%%%%%%%%%%%%%

\end{document}